\documentclass{amsart}

\usepackage{amssymb,mathrsfs,mathdots}
\usepackage[all]{xy}
\usepackage{hyperref}
\usepackage{verbatim}
\usepackage{colortbl}
\usepackage{graphicx,subfigure,psfrag}
\usepackage{enumitem}
\usepackage{combelow}

\newtheorem{theorem}{Theorem}[section]

{\theoremstyle{definition}

\newtheorem{example}[theorem]{Example}
}

\newcommand{\R}{\mathbb{R}}

\begin{document}
	
\title[Bifurcations of polynomial functions]{Bifurcations of polynomial functions with diffeomorphic fibers}
	
\author[F. Braun \MakeLowercase{and} F. Fernandes]
{Francisco Braun$^*$ \MakeLowercase{and} Filipe Fernandes$^\dagger$}

\address{$^*$ Departamento de Matem\'{a}tica, Universidade Federal de S\~ao Carlos, 13565-905 S\~ao Carlos, S\~ao Paulo, Brasil} 
\email{franciscobraun@ufscar.br}

\address{$^{\dagger}$ Universidade do Distrito Federal Professor Jorge Amaury Maia Nunes, 70635--815 Bras\'ilia, Distrito Federal, Brasil}
\email{filipematonb@gmail.com}



\thanks{$^*$ Financed, in part, by the S\~ao Paulo Research Foundation (FAPESP), Brasil, grants 2019/07316-0, 2020/14498-4 and 2023/00376-2; and by the National Council for Scientific and Technological Development (CNPq), Brasil, grants 308112/2023-7 and 403959/2023-3.}

\date{\today}

\begin{abstract}
The phenomena that cause a value of a polynomial function to be a bifurcation one are yet to be described when the fibers have dimension higher than $1$. 
In this note, the main result is the construction of a polynomial submersion function of $\R^3$ with connected fibers having a bifurcation value such that close enough to it the fibers are mutually diffeomorphic. 
We also present an example of a polynomial submersion function of $\R^2$ having a bifurcation value such that close enough to it the fibers are mutually diffeomorphic. 
\end{abstract}
	
\maketitle

\section{Introduction}
It recently came to our knowledge a question raised by Joi\cb{t}a and Tib\u{a}r \cite{JT} whether a polynomial function of $\R^3$ must be a locally trivial fibration at a given regular value whose nearby fibers are connected and diffeomorphic. 
Here we negatively answer this question by presenting an explicit example. 

To be more precise, let $f: \R^m \to \R^n$ be a polynomial map. 
We say that $f$ is a \emph{locally trivial fibration at $\lambda \in \R^n$} if there exist a neighborhood $U \subset \R^n$ of $\lambda$ such that the restriction $f|_{f^{-1}(U)} \to U$ is a \emph{trivial fibration}, that is, the diagram 
$$
\xymatrix{
f^{-1}(\lambda) \times U  \ar[r]^{\varphi} \ar[dr]_{\pi_2} & p^{-1}(U) \ar[d]^f \\
  & U}
$$
is commutative, where $\varphi$ is a homeomorphism and $\pi_2$ is the projection at the second variable. 
A value $\lambda \in \R^n$ is said to be a \emph{bifurcation value} when $f$ is not a locally trivial fibration at $\lambda$. 
We denote by $B(f)$ the set of bifurcation values of $f$. 

It is clear that when $\lambda \notin B(f)$ the fibers $f^{-1}(t)$ are mutually homeomorphic for any $t$ nearby $\lambda$. 
But the converse is far from be true even when $\lambda$ is a regular value in the ``simplest'' case $m = 2$ and $n = 1$, as pointed out by Tib\u{a}r and Zaharia \cite{TZ} by means of an explicit example, although not being a submersion in $\R^2$. 
The submersion assumption shows to be irrelevant, because in Example \ref{456} below we present a submersion $f: \R^2 \to \R$ with a bifurcation value at $0$ and such that all the fibers $f^{-1}(t)$, for $t$ in a neighborhood of $0$, are diffeomorphic, althought not connected. 

The set $B(f)$ is only partially understood. 
When $m = 2$ and $n = 1$ it was completely characterized by Tib\u{a}r and Zaharia in the already mentioned paper \cite{TZ} by means of the phenomena of \emph{vanishing} and \emph{splitting}. 
Joi\cb{t}a and Tib\u{a}r \cite{JT1} extended results of \cite{TZ} for the case $m = n + 1$. 
So for fibers of dimension $1$, the characterization of $B(f)$ is done. 
But when the fibers have dimension greater than $1$ it seems that vanishing and splitting alone do not characterize $B(f)$ anymore. 

In the already mentioned paper \cite{JT}, an explicit example of a submersion map $f: \R^4 \to \R^2$ having connected and mutually diffeomorphic fibers nearby $\lambda$ was contructed. 
In particular, nor vanishing nor spliting can occur here. 

Naturally, one can ask the same for functions. 
This is exactly question Question 3.3 in \cite{JT}: ``if $f: \R^3 \to \R$ have a regular value at $\lambda$, then do constancy of the Euler characteristic plus connectedness of the fibers nearby $\lambda$ guarantee that $\lambda \notin B(f)$?''
In Section \ref{R3} we explicitly construct a non-singular polynomial function $f: \R^3 \to \R$ such that $0 \in B(f)$ even so the fibers $f^{-1}(t)$ are connected and diffeomorphic for all $t$ close enough to $0$. 
So we negatively answer this question. 

We finish the paper by taking advantage of our construction to deliver the above mentioned Example \ref{456}.

\section{The answer is no}\label{R3}
The objective of this section is to present a polynomial function $f: \R^3 \to \R$ such that 
\begin{enumerate}
\item $f$ is a submersion 

\item $0 \in B(f)$ 

\item $f^{-1}(t)$ is connected and diffeomorphic to $\R^2$ minus $3$ points for all $t$ close enough to $0$. 
\end{enumerate}
Let $p: \R \to \R$ be a polynomial function having a local minimum $x_m > 0$ and a local maximum $x_M > 0$ both with value $0$, such that $p(c) \neq 0$ for all $c\notin \{x_m, x_b, x_M\}$, where $x_b$ is in between $x_m$ and $x_M$ and $p(x_b) = 0$. 
In particular, the graphic of $p$ at $(x_b, c_0)$ is transversal to the $x$ axis and $p(0) \neq 0$. 
As an example one can take 
$$
p(x) = (x - 3)^2 (x - 2) (x - 1)^2, 
$$
which has local minimum at $3$ and local maximum at $1$, both with value $0$, and annihilates again only at $2$. 

Consider now $\tilde f: \R^3 \to \R$ given by 
\begin{equation}\label{3d}
\tilde f(x,y,z) = y - p(x), 
\end{equation}
where $p$ is a polynomial function as above. 
This function $\tilde f$ is a submersion with fibers diffeomorphic to $\R^2$. 
Moreover, for every $c$ in a small enough interval $(-\epsilon, \epsilon)$, the fiber $\tilde f^{-1}(c)$ intersects 
$$
\Gamma = \{(x,0,0) \ |\ x \in \R\} 
$$ 
at exactly three points. 

We will work from now on only with $c$'s in $(-\varepsilon, \varepsilon)$, assuming further that $p(0) \notin (-\varepsilon, \varepsilon)$. 
In particular, $\tilde f^{-1}(c) \cap \left(\R^3 \setminus \Gamma \right)$ is diffeomorphic to $\R^2$ minus $3$ points. 
If $c \neq 0$, then the intersections of $\tilde f^{-1}(c)$ with $\Gamma$ are transversal. 
And the intersections of $\tilde f^{-1}(0)$ with $\Gamma$ are tangent at $(x_m,0,0)$ and $(x_M,0,0)$, and transversal at $(x_b,0,0)$. 

Now we consider the diffeomorphism $T: \R^3 \setminus \Gamma \to \R^3 \setminus \Gamma$ given by 
$$
T(x,y,z) = \left(x (y^2 + z^2), y, z \right),\quad T^{-1}(x,y,z) = \left(\frac{x}{y^2 + z^2}, y, z \right). 
$$
In particular, the set $T^{-1}\left(\tilde f^{-1}(c) \cap (\R^3\setminus \Gamma)\right)$ is diffeomorphic to $\R^2$ minus $3$ points. 

Finally we define $f: \R^3 \to \R$ by 
$$
f(x,y,z) = \tilde f \left(x (y^2 + z^2), y, z \right). 
$$
Then $f|_{\R^3\setminus \Gamma} = \tilde f \circ T$. 
In particular, $\nabla f$ is nowhere zero away from $\Gamma$. 
Moreover, since $f_y(x,0,0) = \nabla \tilde f(0)\cdot (0,1,0) = 1$ for all $x$, it follows that $\nabla f$ is nowhere zero in $\R^3$, and hence $f$ is a submersion. 

On the other hand, 
\begin{equation}\label{dfr}
f^{-1}(c) = \left(f^{-1}(c) \cap \Gamma \right) \cup T^{-1} \left(\tilde f^{-1}(c) \cap (\R^3 \setminus \Gamma)\right) = T^{-1} \left(\tilde f^{-1}(c) \cap (\R^3 \setminus \Gamma)\right), 
\end{equation}
because $|f(x,0,0)| > |c|$. 
In particular, the fiber $f^{-1}(c)$ is diffeomorphic to $\R^2$ minus $3$ points. 

Now we prove that $0$ is a bifurcation value of $f$. 
First we prove that $0$ is a bifurcation value of 
$$
g = \tilde f|_{\R^3 \setminus \Gamma}. 
$$ 
If this is not the case, that is, $g$ is a locally trivial fibration at $0$, then there exists $\delta > 0$, a homeomorphism $\varphi: \mathcal{F} \times (-\delta, \delta) \to g^{-1}(-\delta, \delta)$, where $\mathcal{F}$ is $\R^2$ minus three points, and $g \circ \varphi$ is the projection at the second coordinate. 
Then we can select a small neighborhood of $(x_m,0,0)$, and consequently an open set $V$ of $\R^3 \setminus \Gamma \subset g^{-1}(-\delta, \delta)$. 
It is easy to get a contradiction now because $V \cap \tilde f^{-1}(c)$ has no holes when $c>0$ but $V \cap \tilde f^{-1}(0)$ does have a hole and so $\varphi^{-1}(V) \times \left(\mathcal{F} \times \{c\}\right)$ will have a hole for all $c \in (-\delta, \delta)$. 

But \eqref{dfr} gives that 
$$
f^{-1}(-\varepsilon, \varepsilon) = T^{-1} \left(g^{-1} (-\varepsilon, \varepsilon)\right). 
$$
From this and from \eqref{dfr} again we conclude that $0$ is a bifurcation value of $f$ as well. 

\begin{example}\label{456}
Let $\tilde h: \R^2 \to \R$ be defined by 
$$
\tilde h(x,y) = \tilde f(x,y,0), 
$$
with $\tilde f$ given in \eqref{3d}. 
Function $\tilde h$ is a submersion but now the fiber $\tilde h^{-1}(c)$ when intersected with the complement of 
$$
\Gamma = \{(x,0)\ |\ x\in \R\} 
$$
has $4$ connected components. 
We consider 
$$
h(x,y) = \tilde h(x y^2,y), 
$$ 
which is a submersion, as above, whose fiber $h^{-1}(c)$ satisfies 
$$
h^{-1}(c) = J^{-1}\left(\tilde h^{-1}(c) \right), 
$$ 
where $J: \R^2 \setminus \Gamma \to \R^2 \setminus \Gamma$ is the diffeomorphism given by 
$$
J(x,y) = \left(x y^2, y\right). 
$$
In particular, the fiber $h^{-1}(c)$ has $4$ connected components for all $c$ nearby $0$. 
That $0$ is a bifurcation value follows easily because $h$ has both vanishing and splitting at $0$. 
\end{example}

\end{document}